\documentclass[reqno, a4paper, draft]{amsart}
\usepackage{amssymb, amscd, amsfonts,  amsthm}

\usepackage[mathscr]{eucal}
\usepackage{graphicx}
\usepackage[all,poly]{xy} 
 \usepackage[all]{xy}

\newtheorem{theorem}{Theorem}[section]

\newtheorem{satz}[theorem]{}

\theoremstyle{definition}
\newtheorem{dixmier}[theorem]{}

\theoremstyle{definition}

\newcommand{\restrict}{\,{\mathbin{\vert\mkern-0.3mu\grave{}}}\,}

\newcommand{\luk}{\L u\-ka\-s\-ie\-wicz}

\newcommand{\remove}[1]{}

\DeclareMathOperator{\Rn}{{\mathbb R^{\it n}}}

\DeclareMathOperator{\McNn}{\mathcal M([0,1]^{\it n})}

\DeclareMathOperator{\McN}{\mathcal M}
\DeclareMathOperator{\conv}{\rm conv}

\DeclareMathOperator{\Mod}{\rm Mod}
\DeclareMathOperator{\Th}{\rm Th}
\DeclareMathOperator{\I}{[0,1]}
\DeclareMathOperator{\cube}{[0,1]^{\it n}}

\DeclareMathOperator{\relint}{\rm relint}
\DeclareMathOperator{\var}{\rm var}

 \title[The semantics of \L ukasiewicz syntactic consequence]
{The differential semantics of \L ukasiewicz syntactic consequence}

\author{Daniele Mundici}

\address[D. Mundici]{Department of
Mathematics  ``Ulisse Dini'' \\
University of Florence\\
Viale Morgagni 67/A \\
I-50134 Florence \\
Italy}

\email{ mundici@math.unifi.it }

\date{\today}

\begin{document}



\begin{abstract}
The classical condition 
``$\phi$ is a semantic consequence of $\Theta$'' in
infinite-valued propositional \L ukasiewicz logic
\L$_\infty$   is refined  using enriched valuations  that
take into account the effect on $\phi$  of the
 stability of the truth-value of all 
$\theta\in \Theta$  under small perturbations  
(or, measurement errors)
of the models of $\Theta$.  The 
  differential properties of the  
functions represented by $\phi$ and by
all  $\theta\in \Theta$ 
 naturally lead to a new  notion of 
 semantic consequence $\models_\partial$
that  turns out to coincide
with syntactic consequence  $\vdash$. 
\end{abstract}

\maketitle

\hfill {\it to Petr H\'ajek}

\bigskip

\section{Prelude:  semantics for H\'ajek propositional basic logic}
Basic logic  (BL)  was  invented
by  H\'ajek to formalize continuous
t-norms. Certain axioms
satisfied by any such t-norm were singled out in
\cite[2.2.4]{haj};  provability of a formula  $\phi$,
as well as provability of  $\phi$ from a set $\Theta$ of
premises,  were
defined via Modus Ponens, in the usual way,
\cite[2.2.17]{haj}.
BL-algebras,  BL-evaluations of formulas,
and  satisfiability,  
 were then defined in \cite[2.3.3]{haj} and \cite[2.3.8]{haj},
and the following completeness theorem
was proved in \cite[2.3.19]{haj}:

\begin{satz}
\label{BL-completeness}
A formula  $\phi$  is provable
iff every BL-evaluation   satisfies  $\phi$.
\end{satz}

The following strong completeness theorem
directly follows from \cite[2.4.3]{haj}:

\begin{satz}
\label{BL-strong-completeness}
For any 
formula   $\phi$ and 
set $\Theta$  of formulas,
 $\phi$  is provable from $\Theta$
 iff every BL-evaluation   satisfying all
$\theta\in \Theta$  also satisfies $\phi$, in symbols,
$\Theta\models_{BL} \phi$.
\end{satz}

Yet in  \cite[2.3.23]{haj} H\'ajek  champions
a different semantics for BL.
Let us agree to say that
$\phi$  is a {\it t-tautology}  if 
$\phi$ is satisfied by every
evaluation   of $\phi$  into a BL-algebra  arising from
a t-norm.   The resulting
t-tautology semantics  is more adherent to 
 the original motivation of BL-logic: for,  H\'ajek's  BL-axioms
in  \cite[Definition 2.2.4]{haj}  are the result of
his contemplation of continuous t-norms.
The  question arises:  do the BL-axioms 
prove all t-tautologies?
The problem whether BL is {\it the}
logic of continuous t-norms is again posed in a final section
(\cite[9.4.6]{haj}).

In the same pages \cite[9.4.1]{haj},  it is noted
that the traditional semantic consequence relation
$\models$ in  \L$_\infty$  fails to be strongly complete.
A counterexample 
is given in  \cite[3.2.14]{haj}; stated otherwise,
 $\models$ is not compact,
despite model-sets  $\Mod(\psi)$ of  \L$_\infty$-formulas
$\psi(X_1\dots,X_n)$ are compact
subsets of the unit $n$-cube $\cube$, and
compactness has a pervasive role in MV-algebra
theory,  \cite{cigdotmun},\cite{mun11}.

One is then left with two
 rather similar problems involving the mutual role
 of syntax vs. semantics in BL and in  \L$_\infty$:

\begin{itemize}
\item[(A)] {\it Fixed semantics, amendable axioms.}
In case  BL were not complete for t-tautology semantics,
how to  strengthen the BL-axioms  to obtain  a 
strongly complete logic
for continuous t-norms?
\item[(B)] {\it Fixed axioms, amendable semantics.}
It being ascertained  that
$[0,1]$-va\-l\-u\-a\-tions fail to yield a strongly complete
semantics for  \L$_\infty$, 
what new notion of ``model'' of a set of   \L$_\infty$-formulas,
should be devised to get a  strongly complete semantics?
\end{itemize}

In \cite{haj-soft}
H\'ajek himself gave the first substantial contribution
to Problem (A), by adding to BL two  (admittedly not too simple)
axioms which, at the time of  \cite[2.3.23]{haj} and \cite{haj-soft}
were not  guaranteed to follow from the BL-axioms.
 The redundancy of these two axioms was finally
 proved in \cite[5.2]{cigestgodtor}, thus solving
 Problem (A) in  the best possible way: the logic
 originally invented by H\'ajek  
  is indeed strongly complete for valuations
 in t-algebras, the subset of
  BL-algebras directly given by continuous t-norms.
 
 Since the strong completeness of
 $\I$-valuations has been settled in the negative,
  and the \luk\ axioms
 are here to stay,  in
  order  to solve Problem (B)
 we are left with no other choice 
 but to modify the semantics of \L$_\infty$,
looking for a novel, genuinely semantical
notion of $[0,1]$-valuation. 
This is our aim in this paper.

\section{Tangents, differentials and semantic 
consequence relations in \L$_\infty$}
We refer to \cite{cigdotmun} and \cite{mun11} for notation
and background on
MV-algebras and
infinite-valued \L u\-k\-a\-s\-i\-e\-w\-icz propositional logic \L$_\infty$.
The set $\mathsf{FORM}_n$ of \L$_\infty$-formulas
in the variables $X_1,\dots,X_n$ has the same definition
as its boolean counterpart. The \luk\ connectives
$\odot,\oplus$ of conjunction and disjunction
are definable in terms of negation $\neg$ and
implication $\to$. While in boolean logic
formulas take their values in the set $\{0,1\},$
\L$_\infty$-formulas are evaluated in the unit real interval
$\I.$  Let  $\mathsf{VAL}_n\subseteq \I^{\mathsf{FORM}_n}$
denote the set of valuations (also known
as evaluations, assignments, models, interpretations,
possible worlds,\ldots). 
 The {\it truth-functionality}
property of  \L$_\infty$  yields
the following crucial identification:

\begin{satz}
\label{identification}
The set 
$\mathsf{VAL}_n$ can be identified with the unit $n$-cube 
 $\I^n\subseteq \mathbb R^n$ via the restriction map
 $
 V\in \mathsf{VAL}_n \mapsto v=V\restrict
 \{X_1,\dots,X_n\} \in \I^{ \{X_1,\dots,X_n\}}=\I^n.
 $
 For any  fixed formula
 $\phi\in \mathsf{FORM}_n$,
the  map
$
 V\in \mathsf{VAL}_n  \mapsto V(\phi)\in \I
$
defines the  function
$
\hat\phi\colon \I^n\to \I \mbox{\,\, by\,\, } \hat\phi(v)=V(\phi).
$
 The continuity and piecewise linearity of $\hat\phi$
easily follow by induction on the number of
 connectives in $\phi$. 
 \end{satz}

 \begin{dixmier}
 \label{bolzano}
 Following Bolzano and Tarski
 (see \cite[footnote on page 417]{tar}),  \L$_\infty$
 is now equipped with the relation $\models$ of
 {\it semantic consequence} by stipulating that
 for all $\Theta\subseteq \mathsf{FORM}_n$
 and $\phi \in \mathsf{FORM}_n,\,\,\,$
 $
 \Theta\models\phi\,\,\, \mbox{ iff }\,\,\,
 \forall v\in \I^n,\,\,\, (\hat\theta(v) =1 \mbox{ for all }
 \theta\in \Theta \Rightarrow \hat\phi(v)=1).
 $
\end{dixmier}

Mutatis mutandis, this notion of consequence
is   gratified by a completeness theorem
in classical logic and in many nonclassical logics
having totally disconnected valuation spaces. 
However,  

\begin{satz}
\label{directional}
The space $\I^n$  of valuations in
\L$_\infty$ is connected. For
  every $\phi\in \mathsf{FORM}_n$,
    valuation $v\in \I^n$ and
unit vector $u\in \mathbb R^n$
such that $\conv(v,v+\epsilon u)\in\I^n$
for all  small $\epsilon>0,$
the directional derivative
$
{\partial \hat \phi(v)}/{\partial u}
$
 exists and varies continuously with $u$, once
 $v$ is kept fixed. 
 \end{satz} 

The following simple example involving formulas
of one variable
already shows that the differential properties of
 $\hat\theta$ for all $\theta\in \Theta$ are  
 ignored by the semantic consequence relation
 $\models$ of \ref{bolzano},
  although they have no less semantical content than
 the truth-value  $\hat\theta(v):$

 \begin{satz}  Suppose  $\Theta\subseteq
 \mathsf{FORM}_1$ is satisfied by a unique
 valuation $v\in \I$, and $1>v\in\mathbb Q$. Suppose 
   $\partial\hat\theta(v)/\partial x^+=0$
    for all 
 $\theta\in \Theta$.
Let  $\phi=\phi(X)$  be a formula with    
   $\hat\phi(v)=1$  and 
$\hat\phi(w)<1$  for all  $w>v.$
  Then $\Theta\models\phi,$
 although $\,\,\partial\hat\phi(v)/\partial x^+<\,0.$
 \end{satz}
 
Intuitively, the hypothesis
means that each  $\theta\in\Theta$ is not only true at
 $v$, but is also   true for all $w>v$  sufficiently close to
 $v$; in other words,  $\theta$  is  ``stably''  true at $v,\,$ even if
 the value of $v$ were known up to a certain small error
  (depending on $\theta$).  Although 
$\phi$  misses this  (fault-tolerant)  stability
 property of all $\theta\in \Theta,$
 $\phi$ is a semantic consequence of $\Theta$,
 $\,\,\,\Theta\models\phi$.
It should be noted that $\Theta\nvdash\phi.$
Similarly, when   $n>1$   and
    $\Theta\subseteq  \mathsf{FORM}_n,$ the 
higher-order stability properties common to  all  $\theta\in\Theta$
may be missing in some semantic consequence  $\phi$ of $\Theta$.
And again, $\Theta\nvdash\phi.$

While
 directional derivatives make no sense in boolean
 logic, by \ref{directional} 
 they  do make sense in \L$_\infty$. Accordingly, in
 \ref{stably-satisfies} 
  we will give a precise definition of
 ``stable'' consequence relation $\models_\partial$ which is
 sensitive to  all higher order differentiability properties
 of formulas and their associated piecewise linear functions.
In Section \ref{section:appendix} this will be generalized
to   arbitrary (possibly uncountable) sets $\Theta$  of formulas.
%
 In \ref{strong-completeness} we prove that 
  \L$_\infty$ is ``strongly complete'' with respect to
 $\models_\partial\,$:  indeed,  $\Theta\models_\partial\phi$  coincides
 with the syntactical consequence relation $\Theta\vdash\phi.$

We then focus on the relative status of 
  $\models_\partial$ with respect to  $\models$.
As noted in \cite[p.100 and
4.6.6]{cigdotmun}, from Chang completeness theorem we have
  
\begin{satz}  
   \label{when-the-two-coincide}
The two sets
$\Theta^{\models}$ and
$\Theta^{\vdash}$ of 
semantic and syntactic consequences
of a set $\Theta$ of formulas coincide   iff
the Lindenbaum algebra  $\mathsf{LIND}(\Theta)$  is
semisimple.
\end{satz}

\begin{dixmier}
\label{sss}
Following 
  Dubuc and Poveda \cite{dubpov},
we say that an MV-algebra   is
{\it strongly semisimple}  if all its principal
quotients are semisimple. 
\end{dixmier} 

Let $\Theta\subseteq  \mathsf{FORM}_n.$
Building on  \cite{busmun-arxiv},
in  \ref{logical-strong-semisimplicity-due}   we 
 observe that
 $\mathsf{LIND}_\Theta$ is strongly semisimple
 iff  $(\Theta\cup\{\psi\})^{\models}=
(\Theta\cup\{\psi\})^{\models_\partial}$ 
for all $\psi\in \mathsf{FORM}_n$.
Further,  
when $\Theta\subseteq \mathsf{FORM}_1$,
  $\mathsf{LIND}_\Theta$ is strongly semisimple
  iff it is semisimple.
Now suppose
$\mathsf{LIND}_\Theta$  is
semisimple, with
 $\Theta\subseteq \mathsf{FORM}_2$.
  Then  $\mathsf{LIND}_\Theta$  is
strongly semisimple iff 
the set 
$\Mod(\Theta)\subseteq \I^{\{X_1, X_2\}}=\I^2$
of valuations satisfying $\Theta$
has no Bouligand-Severi   \cite{bou, sev}
outgoing  rational tangent  vector at any rational point
$v\in \Mod(\Theta)$.  
See \ref{two-dimensional}.
As shown in \ref{tangent-implies-failure}, 
the existence of a Bouligand-Severi  rational  outgoing tangent
at some rational point $v$ of $\Mod(\Theta)$
entails failure of
strong semisimplicity in the semisimple
MV-algebra   $\mathsf{LIND}(\Th(\Mod(\Theta)))$.

In a final section
 Problems (A) and (B) are retrospectively considered
 in the light of the results of the previous sections. 


 \section{Semantic consequence $\models$ and stable consequence $\models_\partial$} 
 \label{section:consequence}
The  following 
corollary
of Chang's completeness theorem 
is proved in  \cite[3.1.4]{cigdotmun}: 

\begin{satz}
\label{birkhoff}  For each
$n=1,2,\dots,$ the free
$n$-generator MV-algebra  $\McNn$  consists of
all functions  $f\colon \I^n\to \I$ that are 
obtainable from
the coordinate functions  $\pi_i(x_1,\dots,x_n)=x_i$
by pointwise application of the MV-algebraic operations
of negation $\neg x=1-x$
  and truncated addition $x\oplus y=\min(1,x+y).$
  As already noted in \ref{identification},
  any such function $f$ is continuous and piecewise linear.
\end{satz}

  For any
  nonempty closed set $X\subseteq [0,1]^n$
 we let   $\mathcal M(X)$  denote the MV-algebra of
 restrictions to $X$ of the functions in $\McNn$, in symbols,
 $\mathcal M(X)=\{f\restrict X\mid f\in \McNn\}.$
   McNaughton's characterization   \cite[9.1.5]{cigdotmun} 
   of  the free MV-algebra $\McNn$
  will find no use in this paper.

In \cite[3.6.7]{cigdotmun} one can find a proof of 
the following result, which follows
from the proof of Chang's completeness theorem:

 \begin{satz}
 \label{most-general-semisimple}
$\mathcal M(X)$ is a semisimple
 MV-algebra---actually, up to isomorphism, $\mathcal M(X)$
 is the most general possible $n$-generator semisimple
 MV-algebra.
 \end{satz}


For every subset $Y$ of $\cube$,  $\conv(Y)$
denotes the convex hull of $Y$.
To solve Problem (B) we modify  the classical
notion of  valuation as follows:

\begin{dixmier}
\label{differential-valuation}
For  $n=1,2,\dots$ and $0\leq t \leq n$
let  $U = (u_{0}, u_{1},
\ldots , u_{t})$ 
be a  $(t+1)$-tuple of elements of ${\mathbb R}^{n}$
where
$u_{1}, \ldots , u_{t}$ are 
pairwise orthogonal unit vectors.
For  each  $m=1,2,\dots$ let   the
$t$-simplex $T_{U,m}\subseteq \mathbb R^n$ is 
defined  by
\begin{equation}
\label{equation:simplex}
T_{U,m}=
\conv(u_{0},
           u_{0}+u_{1}/m,\,\,
           u_{0}+u_{1}/m+u_{2}/m^2,
           \ldots, u_{0}+u_{1}/m+\cdots +u_{t}/m^t).
\end{equation}
We say that $U$  is a
 {\it differential valuation
  (of order $t$, in $\mathbb R^n$)}
if for all large $m$  the $n$-cube
$\I^n$  contains $T_{U,m}.$
When this is the case,
the set $\mathfrak p_{U} \subseteq \McNn$
is defined by
$ \mathfrak p_{U}=\{f\in  \McNn\mid
f^{-1}(0) \supseteq T_{U,m} \mbox{ for some $m$} \}.$
%
%
\end{dixmier}

Traditional valuations coincide with
differential valuations of order 0.
%



\begin{satz}
\label{prime-ideal-characterization}
Let $U =(u_{0}, u_{1},
\ldots , u_{t})$ be a differential valuation
in $\mathbb R^n$. 

\begin{itemize}
\item[(i)] For all $m=1,2,\dots,$
$T_{U,m}\supseteq T_{U,m+1}$.

\smallskip
\item[(ii)] For every
$\epsilon_{1},\ldots,\epsilon_{t} >0$
there is  $m=1,2,\ldots$  such that
the simplex
$$S = {\rm conv}\{u_{0},
u_{0}+\epsilon_{1}u_{1},\,\,
u_{0}+\epsilon_{1}u_{1}+\epsilon_{2}u_{2},
 \,\ldots, \,\, u_{0}+\epsilon_{1}u_{1}+
 \cdots +\epsilon_{t}u_{t}\}$$
 contains $T_{U,m}$.

\smallskip
\item[(iii)]
$\mathfrak p_{U}$ is a prime ideal of $ \McNn$.

\smallskip
\item[(iv)]
Every prime ideal $\mathfrak p$ of $ \McNn$
has the form $\mathfrak p = \mathfrak p_{V}$ for some
differential valuation $V$.
\end{itemize}
\end{satz}

\begin{proof}  (i)-(ii) are easily verified by induction.
For (iii)-(iv) use (ii) and see 
 \cite[2.8, 2.18]{busmun-apal}. 
\end{proof}

For every convex set $E\subseteq \cube$  we
let $\relint(E)$  denote its relative interior.
The prime ideals  $\mathfrak p_{U}$ of $\McNn$
are conveniently visualized as follows:

\begin{satz}
\label{one-two-three-t}
Let $\,\,U =(u_{0}, u_{1}, \ldots , u_{t})$ 
be a differential valuation in $\mathbb R^n$.  
We then have:

\begin{itemize}
\item[(0)] $\mathfrak p_{(u_0)}$ is the maximal
ideal of  $\McN(\I^n)$ given by all
  functions of  $\McN(\I^n)$ that vanish at $u_0.$

\smallskip
\item[(1)] $\mathfrak p_{(u_0, u_1)}$
is the prime
ideal of  $\McNn$ given by all functions $f\in \McN(\I^n)$
vanishing on  an interval 
  of the form  $\conv(u_0, u_0+u_1/m)$
for some integer $m>0.$ Equivalently,
$f(u_0)=0$ and
$\,\,{\partial f(u_0)}/{\partial u_1}=0.$

\smallskip
\item[(2)]  $\mathfrak p_{(u_0,u_1,u_2)}$  is the prime
ideal  of  $\McNn$  given by those 
$f\in \McN(\I^n)$ 
such that for some integer
  $m>0,\,\,\,$  
  $f$ vanishes on 
  the segment $\conv(u_0,u_0+u_1/m)$, 
  and  $\partial f(y)/\partial u_2=0$ 
for all $y\in \relint\conv(u_0,u_0+u_1/m).$



\medskip
\item[]   And inductively,

\medskip
\item[(t)]
$\mathfrak p_{(u_0,u_1,\dots, u_t)}$  is the prime
ideal of   $ \McNn$ consisting of all 
$f  \in \McN(\I^n)$ such that for some integer  $m>0,\,\,\,$
$f$  vanishes on
 the  $(t-1)$-simplex 
$$S=
\conv\left(u_0,u_0+u_1/m, u_0+u_1/m+u_2/m^2,\ldots, u_0+u_1/m+
\cdots+u_{t-1}/m^{t-1}\right),$$
 and  $\,\,{\partial f(y)}/{\partial u_{t}}=0$ 
  for all 
  $y\in 
  \relint(S).$   
\end{itemize}
\end{satz}

\bigskip
Observe that $\mathfrak p_{(u_0)} \supseteq
\mathfrak p_{(u_0,u_1)}\supseteq\cdots \supseteq
 \mathfrak p_{(u_0,u_1,\dots, u_{t-1})}
\supseteq 
 \mathfrak p_{(u_0,u_1,\dots, u_t)}.$
 
 \bigskip
 Generalizing the classical definitions we can now write:
 
\begin{dixmier}
\label{subbasis-satisfies}
Let $U =(u_{0}, u_{1},
\ldots , u_{t})$ be a differential valuation in $\mathbb R^n$. 
Let $\psi(X_1,\dots,X_n)$ be a formula.
We then say that $U$  {\it    satisfies} $\psi$
if   $1-\hat\psi\in \mathfrak p_U$. Thus  
$$
\hat\psi(u_0)=1,\,\,\, \frac{\partial \hat\psi(u_0)}{\partial u_i}=0,\,\,\,
\dots,
\mbox{ and  $\hat\psi$  satisfies Conditions (2) through (t) 
in \ref{one-two-three-t}}.
$$ 
 \end{dixmier}
 
%
%

\medskip

 \begin{dixmier}
 \label{stably-satisfies}
 For $\Theta\subseteq   \mathsf{FORM}_n$
and $\psi \in   \mathsf{FORM}_n$
we  say that
$\psi$  is a {\it stable consequence} of $\Theta$  and we write
$$
\Theta \models_\partial \psi
$$
if $\psi$ is satisfied by every differential valuation
 $(u_0,u_1,\dots, u_t)$   that    satisfies every $\theta\in \Theta$. 
 \end{dixmier}
 
  \medskip
 
 Observe that 
 $\Theta\models \psi$ in the sense of
 \ref{bolzano}
iff  $\psi$  is satisfied by
every differential valuation
of order 0
satisfying $\Theta$.
Therefore,
\begin{satz}
\label{stable-implies-semantic}
Let   $\Theta\subseteq   \mathsf{FORM}_n$
and $\psi \in   \mathsf{FORM}_n.$
If $\Theta\models_\partial \psi$  then
$\Theta\models \psi$.
\end{satz}


The strong completeness property of 
the stable consequence relation  $\models_\partial $
amounts to the following:

\begin{satz}
\label{strong-completeness}
 $\Theta\models_\partial \psi\,\,\,$  iff
$\,\,\,\Theta\vdash \psi$.
\end{satz}

\begin{proof}
Following \cite[4.2.7]{cigdotmun}, let  $\mathfrak j_\Theta=
\langle \{1-\hat\theta \mid \theta\in \Theta\}\rangle $ be the ideal 
of $\McNn$  generated by the functions given  by all 
 negations of
formulas in $\Theta$.
Equivalently,  $\mathfrak j_\Theta$ is the ideal
generated by the congruence  $\equiv_\Theta$  of
\cite[1.11]{mun11}. 
 Then

\begin{eqnarray*}
\Theta\vdash \psi &\Leftrightarrow & 1-\hat\psi \in \mathfrak j_\Theta,
\mbox{ \cite[4.2.9]{cigdotmun}  or  \cite[1.9]{mun11} }\\
&\Leftrightarrow & 1-\hat \psi \mbox{ belongs to every prime
ideal  $\mathfrak p\supseteq \mathfrak j_\Theta$,  \, 
by subdirect}\\
& & \mbox{ representation, \cite[1.2.14]{cigdotmun}  }\\
&\Leftrightarrow & 1-\hat \psi  \mbox{ belongs to 
every prime  $\mathfrak p$  such that }
1-\hat\theta \in \mathfrak p  \mbox{ for all }  \theta\in \Theta,\\
& & \mbox{ by definition of $\mathfrak j_\Theta$}\\
&\Leftrightarrow &
\mbox{for every differential valuation $U$
in $\mathbb R^n$, if $1-\hat\theta \in \mathfrak p_U$
for all $\theta\in \Theta$}\\
& &\mbox{  then $ 1-\hat \psi \in \mathfrak p_U$,
by \,\ref{prime-ideal-characterization} (iii)--(iv)}\\
&\Leftrightarrow & \psi  \mbox{ is satisfied by all differential valuations }
 U  \mbox{  satisfying all } \theta\in \Theta,  \, \\
 & & \,\,\mbox{by  \ref{subbasis-satisfies}}\\
 &\Leftrightarrow & \Theta\models_\partial \psi,   \,\,\,\mbox{i.e., $\psi$
 is a stable consequence of $\Theta,$
 by  \ref{stably-satisfies}}. 
\end{eqnarray*}
 \end{proof}

\medskip

The ``finitary'' character 
of $\models_\partial$, as opposed to the non-compactness of
$\models,\,\,$  is made precise by the following corollary of
\ref{strong-completeness}:

\medskip
  
 \begin{satz}
\label{finiteness}
Let   $\Theta\subseteq   \mathsf{FORM}_n$
and $\psi \in   \mathsf{FORM}_n.$ Then
 $\Theta\models_\partial \psi\,\,$  iff 
$\,\,\{\theta_1,\dots,\theta_k\}\models_\partial \psi$
for some finite
 subset  $\{\theta_1,\dots,\theta_k\}$ of $\,\Theta$.
\end{satz}

\medskip
Since  $\mathsf{FORM}_n\subseteq \mathsf{FORM}_{n+1}$,
one might ask if $\Theta\models_\partial \psi$ depends
on $n$, so that a more accurate notation would be
 $\Theta\models_{(n,\partial)} \psi$.
  The following immediate 
corollary of \ref{strong-completeness} shows that
such extra notation is unnecessary:
 \begin{satz}
\label{independence}
Let   $\Theta\subseteq   \mathsf{FORM}_n$
and $\psi \in   \mathsf{FORM}_n.$  Then for
any  $m\geq n$,
 $\Theta\models_{(n,\partial)} \psi$
 iff  $\Theta\models_{(m,\partial)} \psi$.
\end{satz}

 \section{Strong semisimplicity
 and  $\models_\partial$} 
 \label{section:tangent-logic}
 Recall from  \ref{sss}  the definition of
strongly semisimple MV-algebra.
Since  $\{0\}$ is
a principal ideal of $A$,  every 
strongly semisimple MV-algebra is  
semisimple.

\begin{satz}
All boolean algebras are strongly semisimple, and so are
 all simple and all finite MV-algebras.
 \end{satz}
 \begin{proof}
 Boolean algebras are hyperarchimedean \cite[6.3]{cigdotmun}.
 The second statement follows from
  \cite[3.5 and 3.6.5]{cigdotmun}. 
 \end{proof}

 \medskip
  The set $\Theta^{\models_\partial}\subseteq 
  \mathsf{FORM}_n$ is  defined by
$
 \Theta^{\models_\partial}=\{\psi\in   \mathsf{FORM}_n
 \mid \Theta \models_\partial \psi\}.
$
 
 \medskip
 
 \begin{satz}
\label{logical-strong-semisimplicity-uno}
Let $\Theta\subseteq \mathsf{FORM}_n$.
Then
$\mathsf{LIND}(\Theta)$  is  semisimple
iff 
$\Theta^{\models} =
\Theta^{\models_\partial}=\Theta^\vdash$.
Thus $\mathsf{LIND}(\Theta)$  is  not semisimple 
iff there is $\psi\in  \mathsf{FORM}_n$ such that
 every differential valuation {\rm of order $0$} 
satisfying  $\Theta$ satisfies $\psi$,  and there is 
a differential valuation 
$U$ satisfying
$\Theta$  but not  $\psi$.
\end{satz}

\begin{proof}
\cite[p.100]{cigdotmun} and \ref{strong-completeness}
 above.
 \end{proof}

 \begin{satz}
\label{logical-strong-semisimplicity-due}
Let $\Theta\subseteq \mathsf{FORM}_n$. Then
$\mathsf{LIND}(\Theta)$  is strongly semisimple
iff for all $\psi\in \mathsf{FORM}_n,\,\,\,$
$(\Theta\cup\{\psi\})^{\models} =
(\Theta\cup\{\psi\})^{\models_\partial}$.
\end{satz}

\begin{proof} 
   For any MV-algebra $A$ and
ideal $\mathfrak j$ of $A$, the quotient
map  
$$\mathfrak i\mapsto \mathfrak i/\mathfrak j
=\{{b}/{\mathfrak i}\mid b \in \mathfrak i\}$$
determines a  1-1 correspondence
 between
 ideals of $A$ containing $\mathfrak j$
 and
  ideals  of
$A/\mathfrak j$,  \cite[1.2.10]{cigdotmun}.
A well known result in universal algebra, 
\cite[3.11]{coh}, 
yields an isomorphism
\begin{equation}
\label{equation:III}
\frac{a}{\mathfrak i}\in \frac{A}{\mathfrak i} \mapsto
 \frac{a/\mathfrak j}{\mathfrak i/\mathfrak j}\in
 \frac{A/\mathfrak j}{\mathfrak i/\mathfrak j}.
\end{equation}
For any   $S\subseteq A$ let $\langle S\rangle$
denote the (possibly not proper) ideal of $A$
generated by $S$.  When $S$ is a singleton
$\{a\}$ we write  $\langle a \rangle$ instead of
$\langle \{a\} \rangle$. 
For  $\mathfrak j$  an ideal
of $A$ we  use the self-explanatory notation
$S/\mathfrak j$  for  $\{b/\mathfrak j\mid b\in S\}$.
For any $a\in A$  we have the trivial identity
\begin{equation}
\label{equation:augmentation}
\frac{\langle a \rangle}{\mathfrak j} = 
\frac{\langle\mathfrak j\cup\{a\}\rangle}{\mathfrak j}.
\end{equation}
For any element  $a/\mathfrak j\in A/\mathfrak i$,
letting  $\langle a/\mathfrak j\rangle$ be the 
 ideal 
generated in $A/\mathfrak j$  by   $a/\mathfrak j$,
a  routine exercise shows  
\begin{equation}
\label{equation:commutation}
\left\langle{a}/{\mathfrak j}\right\rangle=
{\langle a \rangle}/{\mathfrak j}=
\left\{{b}/{\mathfrak j} \mid b \leq m\centerdot a \mbox{ for some }  m=
1,2,\dots \right\}.
\end{equation}
Here are using   the notation
$m\centerdot a$
of   \cite[p.33]{cigdotmun} or
\cite[p.21]{mun11} 
 for $m$-fold truncated addition. 
 
 To complete the proof, for any
$\Theta'$ with $\Theta\subseteq \Theta'\subseteq
\Theta^\vdash$ we have 
$\mathsf{LIND}(\Theta) = \mathsf{LIND}(\Theta')=
\mathsf{LIND}(\Theta^\vdash), $ whence it is no
loss of generality to assume  $\Theta=\Theta^\vdash.$
The set $\{1-\hat\theta\mid \theta\in \Theta\}$ is
automatically an ideal  $\mathfrak j_\Theta$ of $\McNn$ and
we have the isomorphism  
$$
\iota\colon
 \frac{\psi}{\equiv_\Theta} 
\in
\mathsf{LIND}(\Theta)\,\,\, \cong\,\,\,
\frac{1-\hat\psi}{\mathfrak j_\Theta}
\in
 \frac{\McNn}{\mathfrak j_\Theta}\,.
$$
It follows that
the  principal ideal $\langle {\psi}/{\equiv_\Theta} \rangle$
of
$\mathsf{LIND}(\Theta)$ generated by the element
  ${\psi}/{\equiv_\Theta}\in \mathsf{LIND}(\Theta)$ 
corresponds via $\iota$  to the principal ideal
$\langle{ (1-\hat\psi)}/{\mathfrak j_\Theta}\rangle$
generated by the element
$\iota({\psi}/{\equiv_\Theta})=(1-\hat\psi)/\mathfrak j_\Theta
\in \McNn/\mathfrak j_\Theta$.
By (\ref{equation:augmentation})-(\ref{equation:commutation})
 we have the identities
$$
\left\langle\frac{ 1-\hat\psi}{\mathfrak j_\Theta}\right\rangle =
\frac{\langle1-\hat\psi\rangle}{\mathfrak j_\Theta}=
\frac{\langle \mathfrak j_\Theta
\cup \{1-\hat\psi\}\rangle}{\mathfrak j_\Theta}.
$$

   \smallskip
   \noindent
Therefore, 
$\mathsf{LIND}(\Theta)$
is strongly semisimple\,\,\, iff\,\,\, so is
${\McNn}/{\mathfrak j_\Theta}$\,\,\, iff\,\,\, for any principal ideal
${\langle \mathfrak j_\Theta\cup \{1-\hat\psi\}\rangle}/{\mathfrak j_\Theta}$
  of  ${\McNn}/{\mathfrak j_\Theta},$
   the quotient

$$
\frac{{\McNn}/{\mathfrak j_\Theta}}{{\langle \mathfrak j_\Theta\cup \{1-\hat\psi\}\rangle}/{\mathfrak j_\Theta}}\cong
\frac{\McNn}{\langle \mathfrak j_\Theta\cup \{1-\hat\psi\}\rangle}
$$

\medskip 
\noindent
is semisimple.  We are using (\ref{equation:III}).
This is the same as saying  that
$\mathsf{LIND(\Theta\cup\{\psi\})}$ is semisimple
for every  $\psi\in \mathsf{FORM}_n.$
Now apply
\ref{logical-strong-semisimplicity-uno}.
\end{proof}

%

\begin{satz}
For every finite set of \L$_\infty$-formulas  $\Phi$,
the Lindenbaum algebra  $\mathsf{LIND}_\Phi$ is strongly semisimple.
\end{satz}
\begin{proof}
In view of \ref{logical-strong-semisimplicity-due},
 this is  a reformulation of
a  result by Hay  \cite{hay} and W\' ojcicki \cite{woj}
(also see \cite[4.6.7]{cigdotmun} and \cite[1.6]{mun11}), 
stating  that every  
 finitely presented MV-algebra   is  strongly semisimple.
\end{proof}

\medskip
By a quirk of fate,  when
$n=1$
strong semisimplicity boils down to semisimplicity
(see   {\rm (\cite{busmun-arxiv})} for a proof):

\begin{satz}
\label{one-dimensional}
Let $\Theta\subseteq \mathsf{FORM}_1.$
Then 
$\mathsf{LIND}(\Theta)$  is strongly semisimple
iff it is semisimple.
\end{satz}

\medskip
\section{Strong semisimplicity, $\models_\partial$ and
Bouligand-Severi tangents}
While the strong semisimplicity of
$\mathsf{LIND}(\Theta)$  is formulated in purely
algebraic terms,  a deeper understanding of this property
follows from  an exploration  of the 
tangent space of $\Mod(\Theta)$  as a
compact subset of euclidean space  $\mathbb R^n.$

A point $x\in \mathbb R^n$  is said to be
 {\it rational}  if so are all its coordinates.
 By a {\it rational vector}   we mean
 a nonzero vector $w\in \Rn$ such that the line
 $\mathbb Rw=\{\lambda
 w\in \Rn \mid \lambda \in \mathbb R\}\subseteq \Rn$ contains 
 a rational point of $\Rn$ other than the origin.
 Any nonzero scalar multiple of a rational vector is a rational vector.
 
 As usual,  $||v||$ is the length of vector  $v\in \Rn$.
 
 \smallskip
 The following definitions go back to the late twenties and early thirties
 of the past century, and prove very useful to understand the
 geometry of strong semisimplicity, and its relationship with
 stable consequence: 

\begin{dixmier}
\label{definition:halfline-tangent}
{\rm (\cite[\S 53, p.59 and p.392]{sev0},
 \cite[\S1, p.99]{sev},
 \cite[p.32]{bou})}  A half-line $H\subseteq \mathbb R^n$
 is  {\it tangent} to a
 set $X\subseteq \mathbb R^n$ at an accumulation point $x$
 of $X$ if for all $\epsilon, \delta>0$ there is $y\in X$ other than $x$
such that  $||y-x||<\epsilon$, and the angle between $H$ and
 the half-line through $y$ originating at $x$ is $<\delta$.
 \end{dixmier}

\begin{dixmier}
\label{definition:vector-tangent}
   {\rm (\cite[p.16]{book})}
Let  $x$ be an element of   a  closed subset $X$ of $\Rn,$  
and $u$  a unit vector in $\mathbb R^n.$ 
We then say that $u$ is a  {\it 
 Bouligand-Severi  tangent  (unit) vector to $X$    at   $x$}
 if $X$ contains  a sequence $x_0,x_1,\dots$ of elements,
 all different from $x$,  such that
$$
\lim_{i\to \infty }x_i =  x\,\,\, \mbox{ and } \,\,\,\lim_{i\to \infty } {(x_i-x)}/{||x_i-x||}= u.
$$
 \end{dixmier}
 

\noindent
   We further say that $u$  is   {\it outgoing }  if  the open interval   $\relint(\conv(x,x+\lambda u))$
is disjoint from $X$  for some
 $\lambda >0$.

\begin{satz}
 {\rm  (\cite[\S 5, p.103]{sev}).}
 For any  nonempty closed subset $X$ of $\Rn$,
point  $x\in X$,  and   unit vector  $u\in \Rn$
 the following conditions are equivalent:
 \begin{itemize}
\item[(i)]  For all $m=1,2,\dots$,  the cone
\begin{equation}
\label{equation:cone}
 \mathcal C_{x,u,1/m, 1/m^2}
  \end{equation}
with  apex $x$, axis
 parallel to $u$, height $1/m$
and  vertex angle $1/m^2$   
contains infinitely many points of $X$. 

\smallskip
 \item[(ii)]   $u$ is a Bouligand-Severi tangent vector to $X$ at $x$.
 
 \smallskip
 \item[(iii)]   The half-line $x+\mathbb R_{\geq 0}u$
   is tangent to $X$.
 \end{itemize}
 \end{satz}


\begin{satz}
\label{two-dimensional}
{\rm (\cite{busmun-arxiv})}
Let $\Theta\subseteq \mathsf{FORM}_2.$
Suppose
$\mathsf{LIND}(\Theta)$  is  semisimple.
Then
$\mathsf{LIND}(\Theta)$  is strongly semisimple
iff  $\Mod(\Theta)$ does not have any Bouligand-Severi outgoing
rational tangent vector at any of its rational points. 
\end{satz}

\medskip
Combining \cite{busmun-arxiv} 
with our characterization \ref{logical-strong-semisimplicity-uno}
we get

\begin{satz}
\label{tangent-implies-failure}
Let $\Theta\subseteq \mathsf{FORM}_n.$
Suppose 
$\mathsf{LIND}(\Theta)$  is  semisimple and
$\Mod(\Theta)$ has some Bouligand-Severi 
outgoing rational tangent vector $u$ at some rational point
$v\in \Mod(\Theta)$.  Then 
$\mathsf{LIND}(\Theta)$  is not strongly semisimple.
There are formulas  
$\gamma, \lambda \in \mathsf{FORM}_n$  
such that
$\Theta\cup\{\gamma\}\models \lambda$  but it is not the case that
$\Theta\cup\{\gamma\}\models_\partial \lambda$. 
Specifically, while every stable consequence  $\psi$ of
 $\Theta\cup\{\gamma\}$ satisfies
 $\hat\psi(v)=1$  and $\partial \hat\psi(v)/\partial u=0,$
for $\lambda$ we have  $\hat\lambda(v)=1$  and
 $\partial\hat \lambda(v)/\partial u<0.$
 \end{satz}

As in \cite[1.3, 1.4]{mun11},
 the operator $\Th\colon X\subseteq \cube
 \mapsto \Th  X\subseteq \mathsf{FORM}_n$ is
 defined by
 $$
 \Th X =\{\psi\in \mathsf{FORM}_n\mid \hat\psi(w)=1
 \mbox{ for all }  w\in X\}.
 $$

\begin{satz} 
If there exists a Bouligand-Severi rational outgoing tangent
vector at some rational point $v$ of $\Mod(\Theta)$
then  $\mathsf{LIND}(\Th(\Mod(\Theta)))$
is semisimple but  not strongly semisimple.
\end{satz}

\begin{proof} 
 The MV-algebra
 $\mathsf{LIND}(\Th(\Mod(\Theta)))$ is
 semisimple  because  $\Th(\Mod(\Theta))=\Theta^{\models}.$
It is not strongly semisimple by \cite{busmun-arxiv}.
\end{proof}

Thus the strong  semisimplicity of  
$\mathsf{LIND}(\Th(\Mod(\Theta)))$,
and more generally, of every 
$\Phi\subseteq \mathsf{FORM}_n$
with $\Phi^{\models}=\Phi^{\models_\partial},$
only depends on the  (tangent space of the) set
 $\Mod(\Theta)\subseteq \cube$.
 
\section{Concluding remarks}
As shown by the examples of BL and \L$_\infty,\,\,$in the beginning 
we are given a syntactic consequence
relation  $\vdash$ based on a set $\mathcal R$ of 
axioms  and rules.  Then 
variously defined  ``semantic'' consequence 
relations  are tailored around
$\vdash$, until a 
  strongly complete semantic 
consequence relation
is obtained in terms of a certain
set $V^*$  of valuations: 
 in  the case of BL, $V^*$
turns out to be  the subset of
BL-valuations given by  t-algebraic valuations;
in the case of   \L$_\infty$,
$V^*$ is the set of differential valuations,
 which  contains 
 the set of $\I$-valuations
as the special 0-order case. 

Historically, the emergence of semantical notions in 
 first-order logic followed a similar path.
Here a  long distillation process culminated  
in a definitive consequence relation  $\vdash$.
At a later stage, 
 motivation/confirmation of the definitive nature
 of $\vdash$  would be  provided by 
suitably defined  ``models''
    (interpretations, substitutions, evaluations, possible worlds,...). 
  Without them one cannot even speak of the correctness
   of the set $\mathcal R$ of rules of 
first order logic.
The   completeness problem had 
 a long gestation period. 
 The notions of categoricity and completeness
 of theories were often confused with the completeness
 of the set $\mathcal R$ of rules. 
Before the appearance of    Tarskian  models
   over arbitrary universes
the set  of arithmetical models 
  over the {\it fixed} universe
  $\mathbb N$
 was used  to evaluate formulas.

 Turning retrospectively to Problems (A) and (B),
 in the introduction we didn't mention the following
 well known fact  (\cite[20.7]{mun11}):
 
 \begin{satz}
 For each  $i=1,2$  and any (possibly uncountable)
 set $\Theta$ of formulas,  let  
   $\Theta\models_{MVi}\phi$
 be given by the following stipulation:
 \begin{itemize}
 \item[(I)]
    $\Theta\models_{MV1}\phi$
       iff every  $A$-valuation satisfying every $\theta\in \Theta$
   also satisfies $\phi,$  where $A$ ranges over arbitrary
   MV-algebras.
    \item[(II)]
        $\Theta\models_{MV2}\phi$
       iff every  $C$-valuation satisfying every $\theta\in \Theta$
   also satisfies $\phi,$  where $C$ ranges over arbitrary
   MV-chains.
  \end{itemize}
Then   $\models_{MV1}\,\,=\,\,\models_{MV2}\,\,=$ the
syntactic consequence relation $\,\,\vdash$ of 
\L$_\infty$. 
 \end{satz}
 
Each   consequence relation   
 $\models_{MVi}$, while endowing
  \L$_\infty$  with a strongly complete semantics
   has  the same drawbacks as  the
   consequence relation $\models_{BL}$ arising from
   {\it all}   BL-valuations in 
   \ref{BL-completeness}-\ref{BL-strong-completeness}: 
   since  $\models_{MVi}$ does not   
   directly reflect the intuition behind the original
   axioms, its applicability is limited.
   
Consider, for instance, the complexity of the problem  
whether  $\alpha\vdash\beta,$ 
 for $\alpha,\beta\in \bigcup_n \mathsf{FORM}_n.$
 The binary relation  
 $$\vdash_{\rm fin}\,\,\,
 = \,\,\,   \left( \bigcup_n \mathsf{FORM}_n
 \times  \bigcup_n \mathsf{FORM}_n\right)\,\,\bigcap\,\, \vdash$$
 turns out to be  decidable for BL and for \L$_\infty$,  
 no less than for boolean logic.
However, the proper class of {\it all} BL and 
{\it all}  MV algebras,  which is
needed to check
 $\models_{BL}$ and $\models_{MV}$,  
has no role in
the proof of these decidability results.
Actually,  the proof depends on
 subdirect representation
 and completeness theorems, which,
 combined with results like the Hay-W\'ojcicki theorem,
yield a dramatic restriction of the set of  
evaluations  needed to check  semantic consequence.
Suitably small finite chains turn out to be sufficient to decide
if $\beta$  is a consequence of $\alpha.$
In this way
we get   polytime verifiable certificates for 
$\alpha\nvdash \beta$ whence
the coNP-completeness of  
 $\vdash_{\rm fin}$ follows.
 See  \cite{quattro}  and  \cite{mun87}.
 Also see \cite{mon}  for a general discussion of
 strong completeness in various logics, including
  BL and \L$_\infty$.
  
The  evolving semantical notions of
  valuation (model, interpretation, possible world,...),
strongly impinge on the evolution
of  the proof theory of $\vdash.$
While $\vdash$  is immutable,  
the recipe   $\mathcal R$   
to check  $\alpha\vdash \beta$  
is not: we do not even know if
``proofs'', as we understand them today in
boolean logic  (let alone \L$_\infty$  and  BL)
will one day be replaced by revolutionary polytime
decision procedures.
  
%
%

H\'ajek's  intuition of the BL-axioms
 was confirmed by a definitive strong completeness result
 for  valuations over
 t-algebras rather than over arbitrary BL-algebras. 
 Similarly,  the \luk\ axioms for \L$_\infty$, 
 as well as Chang's MV-algebraic
axioms are now gratified by a strongly complete
(genuinely  semantic)
consequence relation $\models_\partial$  
that does  not  resort to valuations over exoteric MV-algebras
and their ``infinitesimal truth-values''.
 Rather, $\Theta \models_\partial\psi$ 
 depends on (real-valued) differential valuations that  
check if $\psi$  has the stability properties common to
all    $\theta\in \Theta$.
 
%
 
Closing a circle of logic-algebraic-geometric ideas,
our results in this paper show that 
 the
 traditional semantic consequence relation $\Theta\models\phi$
 fails to be strongly complete because of its 
  total  insensitivity
 to the Bouligand-Severi  tangent
space of  $\Mod(\Theta)$. 
Strong completeness is retrieved
by differential valuations,  which
take into account the directional
derivatives of formulas along the
 tangent space of $\Mod\Theta$.

\section{Appendix: stable consequence for arbitrary sets of sentences}
\label{section:appendix}
Since MV-algebras are  Lindenbaum algebras
of set of formulas in \L$_\infty$, we have to consider
arbitrarily large sets of formulas on unlimited supplies
of variables. 
So let
 $\mathcal X=\{X_1,X_2,\dots,X_\alpha,\ldots\mid \alpha<\kappa\}$
be a set of variables of infinite,
possibly uncountable  cardinality $\kappa$, indexed by all ordinals $0<\alpha<\kappa.$  We let $\mathsf{FORM}_\mathcal X$ be the
set of formulas $\psi(X_{\alpha_1},\dots,X_{\alpha_t})$ whose variables
are contained in $\mathcal X$.
In this appendix we routinely extend
Definition \ref{stably-satisfies} to arbitrary
  subsets $\Theta$  of $ \mathsf{FORM}_\mathcal X$
and formulas $\psi\in  \mathsf{FORM}_\mathcal X$.

\begin{satz}
The free MV-algebra over $\kappa$  free generators
is the MV-algebra  $\McN(\I^\kappa)$ of all functions
on the Tychonov cube  $\I^{\mathcal X}=\,\I^\kappa$ 
generated by the coordinate functions
$\pi_\beta(x)=x_\beta,\,\,\,(x\in \I^\kappa,\,\,\,
0<\beta<\kappa)$
 by  pointwise application
of the $\neg,\oplus$ operations,  \cite[9.1.5]{cigdotmun}. 
\end{satz}

\begin{dixmier}
\label{first-identification-new}
For any finite set $\mathcal K=\{X_{\alpha_1},\dots, X_{\alpha_d}\}
\subseteq \mathcal X$
we identify $\I^{\mathcal K}$ with  the
 set of all  $x\in \I^\kappa$  such that all coordinates $x_\beta$
of $x$ vanish, 
with the possible exception of  $\beta\in \{\alpha_1,\dots,\alpha_d\}.$
For any formula   $\phi \in \mathsf{FORM}_\mathcal X$,
we let
  $\var(\phi)$ be the set of variables occurring in $\phi$.
Identifying the function $\hat\phi$  with an element of
  $\McN(\I^{\var(\phi)})$ we will tacitly 
identify    $\McN(\I^{\mathcal K})$ 
with the  
subalgebra of $\McN(\I^{\kappa})$ 
consisting of all
$\hat\phi$  for $\phi\in \mathsf{FORM}_{\mathcal K}$.
\end{dixmier}


\begin{dixmier}
\label{order}
Suppose now we are given two finite subsets
$\mathcal H \subseteq \mathcal K \subseteq \mathcal X$ and
two differential valuations  $U=(u_0,u_1,\dots,u_d)$ 
in $\mathbb R^{\mathcal H}$
and $V=(v_0,v_1,\dots,v_e)\in \mathbb R^{\mathcal K}$.
We then have  
two prime ideals  $\mathfrak p_U$
of $\McN(\I^\mathcal H)$ and $\mathfrak p_V$
of $\McN(\I^\mathcal K)\supseteq \McN(\I^\mathcal H).$
Recalling \ref{differential-valuation}, we say that
$V$ {\it dominates }  $U$, in symbols,
$V\succeq U,$  if  $\mathfrak p_U=\mathfrak p_V\cap 
\McN(\I^\mathcal H)$.   Whenever $V\succeq U$, the point  
$u_0$ of  $\I^\mathcal H\subseteq \I^\mathcal K$
is obtained by forgetting all coordinates of $v_0$
other than those in $\mathcal H$.  
Further information on the relationship between
$U$ and   $V$
can be found in \cite[\S 4]{busmun-apal}. 
\end{dixmier}

The following definition is a straightforward generalization
of  \ref{differential-valuation}:

\begin{dixmier}
\label{system}
A {\it differential valuation   in
$\mathbb R^\kappa$}   is a {\it $\succeq$-direct system} 
\begin{equation}
\label{equation:system}
W=\{U_\mathcal H\mid \mathcal H\subseteq \mathcal X,
\,\,\mathcal H \mbox{ finite}\}
\end{equation}
of differential valuations  $U_\mathcal H$ in $\mathbb R^\mathcal H$,
indexed by all finite subsets  $\mathcal H$  of $\mathcal X.$
As usual, directedness means that, for any  finite $\mathcal H',
\mathcal H''\subseteq \mathcal X,\,\,\,$ 
$U_{\mathcal H'\cup\mathcal H''}$ dominates
both $U_{\mathcal H'}$ and
and
$U_{\mathcal H''}.$
We say that
$W$ {\it satisfies} a formula  $\phi\in \mathsf{FORM}_\mathcal X$
if  $U_{\var(\phi)}$ satisfies $\phi$ in the sense of 
\ref{subbasis-satisfies},
i.e.,  $1-\hat\phi$   belongs to $\mathfrak p_{U_{\var(\phi)}}$.
\end{dixmier}

 \begin{dixmier}
 \label{stably-satisfies-bis}
 For $\Theta\subseteq   \mathsf{FORM}_\mathcal X$
and $\psi \in  \mathsf{FORM}_\mathcal X$
we  say that
$\psi$  is a {\it stable consequence} of $\Theta$  and we write
$
\Theta \models_\partial \psi,
$
if $\psi$ is satisfied by every differential valuation
  $W$   in $\mathbb R^{\kappa}$
    that    satisfies every $\theta\in \Theta$. 
 \end{dixmier}
 
Recalling
\ref{independence}, 
it is not hard to see that
$\models_\partial $ is an extension of the stable consequence
relations defined for $\Theta\subseteq \mathsf{FORM}_n$ and
$\psi\in  \mathsf{FORM}_n,\,\,\,(n=1,2,\dots).$ 
The ``strong completeness'' theorem for 
this general consequence relation
$\models_\partial $
now states:

\begin{satz}
\label{strong-completeness-bis}
For any (possibly uncountable) set  $\mathcal X$  of variables,
$\Theta\subseteq \mathsf{FORM}_\mathcal X$
and $\psi\in \mathsf{FORM}_\mathcal X,$
 $\Theta\models_\partial \psi\,$  iff
$\,\Theta\vdash \psi$.
\end{satz}

\begin{proof}
Every prime ideal  $\mathfrak p$  of $\McN(\I^{\kappa})$
is uniquely determined by its intersections
$\mathfrak p\cap \McN(\I^{\mathcal H}) $
letting $\mathcal H$ range over finite subsets
of  $\mathcal X.$  Any such intersection is a
prime ideal of  $ \McN(\I^{\mathcal H}) $.
By  \cite[2.18]{busmun-apal}, for every finite $\mathcal H\subseteq
\mathcal X$
there is a differential valuation
$U_\mathcal H$  in $\mathbb R^\mathcal H$
such that
$\mathfrak p\cap \McN(\I^{\mathcal H})=
\mathfrak p_{U_\mathcal H}$.  Letting  now
$\mathcal H$ range over all finite subsets
of $\mathcal X$,  
 the $\mathfrak p_\mathcal H$ make a
$\supseteq$-direct system with union
$\mathfrak p$.
Correspondingly
the differential valuations
$U_\mathcal H$ in $\mathbb R^\mathcal H$
make a $\succeq$-direct system, i.e., a differential
valuation  $W=W_\mathfrak p$ in $\mathbb R^\kappa$ of the
form  (\ref{equation:system}). 
Every prime ideal $\McN(\I^{\mathcal H})$
arises in this way from a differential valuation
$W_\mathfrak p$  in $\mathbb R^\kappa$.
Now argue as in the proof
of \ref{strong-completeness} using the subdirect
representation theorem for $\McN(\I^{\mathcal H})$.
\end{proof}

\bibliographystyle{plain}

 \end{document}